\input amssym.def
\input amssym
\documentclass{article}

\begin{document}
\bibliographystyle{plain}

\newtheorem{thm}{Theorem}
\newtheorem{defin}[thm]{Definition}
\newtheorem{lemma}[thm]{Lemma}
\newtheorem{propo}[thm]{Proposition}
\newtheorem{cor}[thm]{Corollary}
\newtheorem{conj}[thm]{Conjecture}
\newtheorem{exa}[thm]{Example}

\centerline{\LARGE \bf Measurable Categories}

\vspace*{1cm}

\centerline{\parbox{2.5in}{D.\ N.\ Yetter \\ 
Department of Mathematics \\
Kansas State University \\ Manhattan, KS 66506}}
\vspace*{1cm}

{\small
\noindent{\bf Abstract:} We develop the theory of categories of measurable
fields of Hilbert spaces and bounded fields of operators.  We examine 
classes of functors and natural transformations with good measure
theoretic properties, providing in the end a rigorous construction for
the bicategory used in \cite{CY.M2G} and \cite{CS.2G} as the basis for
a representation theory of (Lie) 2-groups.  Two important technical
results are established along the way:  first it is shown that all invertible 
additive bounded functors (and thus {\em a fortiori} all invertible *-functors)
between categories of measurable
fields of Hilbert spaces are induced by invertible measurable transformations
between the underlying Borel spaces and second we establish
the distributivity of Hilbert
space tensor product over direct integrals over Lusin spaces with respect
to $\sigma$-finite measures. The paper concludes with
a general definition of measurable bicategories.}

\section{Introduction}

One of the vexing problems in the algebraic approach to quantum topology
and related work on algebraic models for quantum gravity has 
been a lack of suitable examples of the algebraic structures.  Just as 
classical examples of tensor categories were ``too commutative'' to have any
but the most trivial relation to classical knots and three manifolds, so 
known examples of algebraic structures of the type expected at a formal level 
to be related to the structure of 4-manifolds appear to be ``too commutative''
to detect anything other than homeomorphism type.

The desire to find sufficiently non-commutative examples of such structures
(for example, monoidal bicategories with appropriate dualities) has been
one motivation for the work in higher dimensional algebra done in the
past decade.  Another has been the suspicion that symmetry groups may be 
inadequate expressions of the symmetries needed to formulate a quantum
theory of gravity.

The present work is intended to address difficulties in one line of
development in this direction:  the representation theory of categorical
groups (or, as they are called when considered as a type of bicategory
with one object, 2-groups\footnote{Not to be confused with the use of 
``2-groups'' in reference to 2-torsion in a finite group}).
This representation theory is not developed in the present paper
beyond stating the definition to give an example of a general
notion of measurable bicategory.  Rather this is
the subject of \cite{CY.M2G} and \cite{CS.2G}.

The difficulties in representing 2-groups with infinite sets of objects in
any of the versions of 2-VECT considered by Kapranov and Voevodsky \cite{KV}
(cf. also \cite{B}, \cite{BM})  
also are analogous to the difficulties in representing non-compact 
groups in the
category of finite-dimensional vector spaces.  The natural way to 
overcome them is the same:  move to a setting built out of measurable spaces
rather than finite sets.

Just as in 2-VECT one has families of vector-spaces indexed over finite sets,
so we need a setting where we have families of Hilbert spaces indexed over
measure spaces.

The appropriate ideas have already been considered in the context of 
functional analysis and the (unitary) representation theory of non-compact
groups.  There the generalization of direct sum decomposition theorems required
the introduction of a measure-theoretic analogue:  the direct integral.

Although the constructions in Sections 2 through 4 are perfectly general
and will work for any measurable space satisfying the 
mild technical hypothesis that points are measurable, 
beginning with section 5, where
we will need to invoke results of Maharam \cite{Maharam} (cf. also
\cite{GrMa}) on disintegrations of measures,
we will require the hypothesis that the measurable spaces in question 
are the Borel space associated to a Lusin space (the image of a separable
metric space under a continuous bijection) whose points are
Borel sets.\footnote{This is equivalent to an odd separation axiom stronger
than $T_0$, but weaker than $T_1$, ``Freyd's Favorite Separation Axiom'': 
every point is locally closed, that is, every point admits an open neighborhood
in which it is a closed set.} Beginning in Section 4, 
although the definitions are applicable to general measures, our results
will require the use of the Radon-Nikodym Theorem and the Lebesgue Decomposition
Theorem.  Thus we will assume throughout that all measures are 
totally $\sigma$-finite.  Finally, we will need at one point to
assume that the $L^2$ spaces for all measure spaces considered
are separable Hilbert spaces.  Observe that this rather long list of technical hypotheses is satisfied by among others the Haar measure on
any finite dimensional Lie group and by any measure on a Lusin space concentrated at a countable set of points.

\section{Categories of Measurable Fields \\ of Hilbert Spaces}

The rich structures associated to the direct integral 
construction appear never to have
been examined from a categorical point of view, although all the 
necessary ingredients are there.  Indeed 
the first two definitions are variations of notions found in Takesaki 
\cite{Takesaki}:

\begin{defin} \label{meas.field}
A {\em measurable field of Hilbert spaces} ${\cal H}$ on a Borel space
$(X,S)$ is a pair $({\cal H}_x, {\cal M}_{\cal H})$, where ${\cal H}_x$ is
an X-indexed family of Hilbert spaces, and ${\cal M}_{\cal H} = {\cal M}$ 
is a linear
subspace of $\prod_{x \in X} {\cal H}_x$ (the product as vector-spaces)
satisfying

\begin{enumerate}
\item $\forall \xi \in {\cal M} \;  x \mapsto \|\xi(x)\|_x$ is measurable,
\item For all $\eta \in \prod_{x \in X} {\cal H}_x$ the measurability of 
	$ x \mapsto \langle \eta(x) | \xi(x) \rangle_x$ for
	all $\xi \in {\cal M}$ implies $\eta \in {\cal M}$,
\item $\exists \{ \xi_i \}_{i = 1}^\infty \subset {\cal M}$ such that
	$\{\xi_i(x)\}_{i = 1}^\infty$ is dense in ${\cal H}_x$ for all 
$x \in X$.
\end{enumerate}

An {\em almost measurable field of Hilbert spaces} ${\cal H}$ on a Borel
space $(X,S)$ is a pair $({\cal H}_x, {\cal M}_{\cal H})$ as above, 
satisfying conditions 1 and 2, but not necessarily condition 3.
\end{defin}

One thing which should be noted immediately is that condition 1, together with 
a polarization argument, shows
that condition 2 is really bidirectional:  $\xi \in {\cal M}$ if and only if 
for all $\zeta \in {\cal M}$, $x \mapsto
\langle \xi(x) | \zeta(x) \rangle_x$ is measurable.

In what follows, we will assume w.l.o.g.
that the sequence $\{\xi_i\}$ verifying condition 3. includes
the 0 section.  This assumption simplifies some constructions.

It might seem logical to use ``measurable field'' to refer to what we have 
called an almost measurable field.  There is, however, reason beyond 
adherence to established terminology not to do so.  Example \ref{pathology} 
shows that although the measurablity conditions are maintained in the 
definition of an almost measurable field, the lack of the countability or 
separability condition allows non-measurable pathologies.

\begin{defin}
A {\em measurable field of bounded operators} $\phi$ from ${\cal H}$ to 
${\cal K}$, for ${\cal H}$ and ${\cal K}$ (almost) measurable fields of 
Hilbert spaces is an
X-indexed family of bounded operators $\phi_x \in B({\cal H}_x, {\cal K}_x)$
such that $\xi \in {\cal M}_{\cal H}$ implies 
$\phi(\xi) \in {\cal M}_{\cal K}$,
where $\phi(\xi)_x = \phi_x(\xi_x)$.

A measurable field of bounded operators is {\em bounded}
if the real valued function $x \mapsto \|\phi_x\|_x$ is bounded.

A measurable field of bounded operators is {\em essentially bounded} with
respect to a measure $\mu$ on $(X,S)$ if
$x \mapsto \|\phi_x\|_x$ is in $L^\infty(X,\mu)$.  (Here $\| \; \|_x$ 
denotes the
operator norm on $B({\cal H}_x, {\cal K}_x)$.)
\end{defin}
Classically measure spaces are considered, and it is thus
more natural to work with essentially bounded fields, as
two field of operators which differ on a set of measure zero 
will induce the same operator between direct integrals.  
We, however, are working in a setting where
different measures will be considered on the same Borel 
space, and are thus obliged to work with
bounded fields, ``measure zero'' having no fixed meaning
when one changes measures.

We can then organize these into a category:

\begin{defin}
The {\em category of measurable fields of Hilbert spaces on} $(X,S)$ has as
objects all measurable fields of Hilbert spaces on $(X,S)$ and as arrows
all bounded fields of operators on $X$.  Source, target, 
identity arrow and composition are obvious.  We denote this category by
${\bf Meas}(X,S)$.

Similarly, the {\em category of almost measurable fields of Hilbert spaces} on
$(X,S)$ has as objects all almost measurable fields of Hilbert spaces on 
$(X,S)$ and as arrows all bounded fields of bounded operator 
between them.  We denote this category by ${\bf AlMeas}(X,S)$.
\end{defin}

It will be important in what follows to organize these families of categories
into 2-categories by introducing suitable functors and natural transformations
between them.

Before we do this, however, we consider an example which explains why
the name ``measurable field'' is a properly applied to the classical notion 
rather than to
the more general notion, and derive some elementary properties of the 
categories themselves.

\begin{exa} \label{pathology}
Consider the almost measurable field of Hilbert spaces on 
$({\Bbb R},S)$, the real line with the Borel structure of all
Borel measurable sets, defined 
as follows:  Let $X$ be a non-measurable subset of 
${\Bbb R}$, and let $\{A_\lambda \}$ be a (necessarily uncountable) set of 
measurable sets such that $\cup_\lambda A_\lambda = X$.  Now, consider the
field of Hilbert spaces 

\[ {\cal H}_x = \left\{ \begin{array}{ll}
				{\Bbb C} & \mbox{if $x \in X$}\\
				0 & \mbox{otherwise}
			\end{array} \right. \]

\noindent with ${\cal M}$ given by $\{ \zeta | \mbox{$x\mapsto 
\langle \zeta | \xi_\lambda \rangle$ is measurable for all $\lambda$} \}$,
where $\xi_\lambda$ is the section which is 1 on each fiber over $A_\lambda$
and 0 otherwise.
\end{exa}

Note the pathology present in the example just given:  although the sections
in $\cal M$ all have measurable fiberwise norms, and give 
measurable functions when their fiberwise scalar products are taken, the 
support of the field is non-measurable.   This is so even though the fibers
are all separable.

\begin{propo}
For any Borel space $(X,S)$, the categories ${\bf Meas}(X,S)$ and 
${\bf AlMeas}(X,S)$ are ${\Bbb C}$-linear additive $C^*$-categories.
\end{propo}

\noindent{\bf proof:} It is easy to see that the hom-set, equipped with
fiberwise addition, and multiplication by scalars are 
modules over the algebra of bounded functions on $X$, 
and thus, {\em a fortiori} vectorspaces over $\Bbb C$. 
Moreover they are equipped with a norm $\|\phi\| = \sup_{x\in X} \|\phi_x\|_x$
where $\| \;\; \|_x$ is the operator norm of the bounded operator at $x$.

The field with constant fiber 0 is plainly a zero object.  It thus
remains to show that the category admits biproducts. 

As observed in the Appendix, {\bf Hilb} is an additive category.  It follows
that ${\bf Hilb}^X$ is as well.  Now observe that there are forgetful
functors from ${\bf Meas}(X,S)$ and ${\bf AlMeas}(X,S)$ to ${\bf Hilb}^X$.  
We claim that these
functors creates biproducts, in the sense that the biproduct 
$\{{\cal H}_x \oplus {\cal K}_x\}$ of the underlying
$X$-indexed families of Hilbert spaces $\{{\cal H}_x\}$ and
$\{{\cal K}_x\}$ has the structure of a(n almost) 
measurable
field of Hilbert spaces, and that the projections and inclusions are
measurable fields of operators.

The fibers of the biproduct are the direct sum ${\cal H}_x \oplus {\cal K}_x$
of the underlying vectorspaces, with the sum of the scalar products on the 
summands as scalar product (see the Appendix).   We can then form
${\cal M}_{{\cal H} \oplus {\cal K}}$ by taking the closure of the set
$G = \{ (\eta, 0) | \eta \in {\cal M}_{\cal H}\} \cup
\{ (0,\kappa) | \kappa \in {\cal M}_{\cal K}\}$
under condition 2. of the definition of (almost) measurable fields
(cf. \cite{Takesaki}).

Now, observe that
if $\{\eta_i\}$ and $\{\kappa_i\}$ are fundamental sequences for $\cal H$
and $\cal K$ respectively, then the sequence $\{(\eta_i, \kappa_j)\}$ 
satisfies the required condition that $\{(\eta_i(x), \kappa_j(x))\}$
is dense in ${\cal H}_x \oplus {\cal K}_x$, and moreover is in
${\cal M}_{{\cal H} \oplus {\cal K}}$. since the scalar products
of elements with elements of $G$ are plainly measurable

Now it is easy to see that the inclusions preserve measurable sections:  for
example, for the first inclusion
$\langle (\zeta, 0) | (\eta, 0) \rangle = \langle \zeta | \eta \rangle$,
while $\langle (\zeta, 0) | (0, \kappa) \rangle = 0$
and thus $(\zeta, 0)$ is measurable whenever $\zeta$ is.

For the projections, consider a section $(\zeta, \lambda)$.  By construction
it is measurable whenever the fiberwise scalar products with each of the
$(\eta, 0)$'s and each of the $(0, \kappa)$'s are measurable functions.
But taking scalar products with these reduces to taking scalar products
of one of the summands alone, thus implying that each of the summands
$\zeta$ and $\lambda$ are measurable.  

For boundedness, it follows from the construction of the biproduct
on the direct sum in ${\bf Hilb}$ that the inclusions (resp. projections) are norm preserving
(resp. norm decreasing) in each fiber, and thus, taking
suprema are norm preserving (resp. norm decreasing). 

It thus remains to show that the hom-sets are complete with respect 
to the norm, that $\|\phi(\psi)\| \leq \|\phi\|\|\psi\|$ and $\|\phi(\phi^*)\| = \|\phi\|$.
All three follow from the corresponding result for ${\bf Hilb}$ once it is
observed that equations and non-strict inequalities are preserved by suprema, and
that Cauchy sequences of bounded fields of operators necessarily give Cauchy sequences
in each fiber since the norm is the supremum of the operator norms in each fiber.
$\bullet$ \smallskip

To further develop the theory, we need some notions
native to categories of the form ${\bf Meas}(X,S)$ and
${\bf AlMeas}(X,S)$.

\begin{defin}
The {\em support} of (an almost) measurable field of Hilbert 
spaces $\cal H$ on a Borel space $(X,S)$ is the set 
$supp({\cal H}) = \{ x \in X | {\cal H}_x \not\cong 0 \}$.

The support of a measurable section $\xi \in {\cal M}_{\cal H}$ 
is the set $supp(\xi) = \{ x \in X | \xi(x) \neq 0 \}$.  Note 
by taking norms that it is necessarily
a measurable set.

The support of a field of bounded operators $B$ is the set
$supp(B) = \{ x \in X | B(x) \neq 0 \}$.
\end{defin}

\begin{defin}
If $\cal H$ (resp. $\xi$, $B$) is a(n almost) measurable field of 
Hilbert spaces (resp. a measurable section, 
a measurable field of operators) on $(X,S)$,
and $A \in S$, then the restriction of $\cal H$ (resp. $\xi$, $B$) to $A$, 
denoted ${\cal H}|_A$ (resp.
$\xi|_A$, $B|_A$), is given by

\[ {\cal H}|_{Ax} = \left\{ \begin{array}{ll}
	{\cal H}_x & \mbox{if $x \in A$} \\
	0 & \mbox{otherwise}
\end{array} \right. \]

\noindent with ${\cal M}_{{\cal H}|_A} = {\cal M}_{\cal H} 
\cap \prod_x {\cal H}|_{Ax}$, where we identify
${\cal H}|_{Ax}$ with a subspace of ${\cal H}_x$, either the entire Hilbert 
space or 0 (resp.

\[ \xi|_A(x) = \left\{ \begin{array}{ll}
	\xi(x) & \mbox{if $x \in A$} \\
	0 & \mbox{otherwise}
\end{array} \right. , \]

\[ B_A(x) = \left \{ \begin{array}{ll}
        B(x) & \mbox{if $x \in A$} \\
	0 & \mbox{otherwise ).}
\end{array} \right. \]
\end{defin}

In the measurable case, 
we need to see that if $\{\xi_i \}$ 
is a fundamental sequence for $\cal H$, 
then $\{\xi_i|A \}$ is a fundamental sequence for 
${\cal H}|A$.  This will follow directly from

\begin{lemma} \label{res.lemma}
For any measurable set $A$ and any (almost) measurable field of Hilbert spaces
$\cal H$, $\xi \in {\cal M}_{\cal H}$ implies $\xi|_A
\in {\cal M}_{\cal H}$.
\end{lemma}

\noindent{\bf proof:} We show that $f(x) = \langle
\xi|_A \; | \zeta \rangle$ is measurable for any $\zeta 
\in {\cal M}_{\cal H}$.  Let $g(x) = \langle
\xi | \zeta \rangle$.  Recall that a real-valued 
function $\phi$ is measurable if for any measurable set 
$M \subset {\Bbb R}$ the set $N(\phi) \cap \phi^{-1}(M)$ is 
measurable, where $N(\phi) = \{x | \phi(x) \neq 0\}$. 
Observing that $f$ and $g$ agree 
on $A$, but that $f$ is constant 0 on $\neg A$, we have
that $N(f) = N(g) \cap A$, while $f^{-1}(M) = g^{-1}(M)
\cap A$ if $0 \not\in M$ and $f^{-1}(M) = g^{-1}(M) \cup \neg A$ if 
$0 \in M$.

In either event the intersection $N(f) \cap f^{-1}(M)$
is $N(g) \cap g^{-1}(M) \cap A$, which is measurable since
both $A$ and $N(g) \cap g^{-1}(M)$ are. $\bullet$ \smallskip

Restrictions in general provide objects ${\bf Meas}(X,S)$ with direct sum 
decompositions:

\begin{thm} \label{res.sum}
If $\cal H$ is a measurable field of Hilbert spaces on
$(X,S)$ and $A \in S$, then

\[ {\cal H} \cong {\cal H}|_A \oplus {\cal H}|_{\neg A} \]

\end{thm}

\noindent{\bf proof:} This follows almost immediately
from three observations about $|_A$ on sections:  first
$|_A$ is idempotent, second $|_A|_{\neg A} = 0$, and
third for any section $\xi$, we have $\xi = \xi|_A + \xi|_{\neg A}$.$\bullet$ \smallskip

\begin{defin}
A measurable field of Hilbert spaces $\cal H$ on $(X,S)$
 is a
{\em partial measurable line bundle} if for all $x \in X$
either ${\cal H}_x \cong {\Bbb C}$ or ${\cal H}_x \cong 0$.
\end{defin}

We have already seen a pathology which can 
occur when the corresponding notion is 
considered in the almost measurable setting.  
We now establish that it does not occur for partial 
measurable line bundles.  Indeed we have

\begin{thm}
$supp(-)$ induces a canonical bijection between the isomorphism classes of 
partial measurable line bundles
on $(X,S)$ and the collection $S$ of measurable sets.
\end{thm}

\noindent{\bf proof:} We proceed in two stages:  
first we show that $supp(-)$ takes values in the measurable sets, and then 
we construct an inverse.

Consider a partial measurable line bundle 
${\cal H}_x, {\cal M}, \{\xi_i\}$.  Now it is 
clear that the union of the supports of any family of measurable sections of 
$\cal H$ is contained in the support of $\cal H$, as
fibers outside the support are 0.  It 
thus suffices to see that there is a 
set of sections whose support contains that of $\cal H$ and is measurable.

Consider the fundamental sequence 
$\{\xi_i \}$.  By the density condition, 
given any non-zero element $v$ of the fiber $\cal H_x$, there is a $\xi_i$ such that 
$\| v - \xi_i(x) \| < \frac{1}{2}\|v\|$, and thus
such that $\xi_i(x) \neq 0$.  From this it 
follows immediately that $\bigcup_i supp(\xi_i) = supp({\cal H})$.

But, the $supp(\xi_i)$'s are a 
countable set of measurable sets in a Borel space, and thus their 
union is measurable.

To construct an inverse to 
$supp(-)$, we proceed as follows:  Given a measurable set $A \subset X$, 
we may form a partial measurable line 
bundle with support $A$ by taking the restriction of the constant measurable 
field $\Bbb C$.

It thus suffices to show that any partial measurable line $\cal H$ bundle with 
support $A$ is isomorphic to ${\Bbb C}|_A$.

To construct the isomorphism, we first construct a measurable section of 
$\cal H$ whose support coincides with that of $\cal H$:

\[ \xi_{total} = \sum_{i = 1}^\infty \xi_i|_{\neg \cup_{j < i} supp(\xi_j)} \]

As usual, to see that this sum lies in $\cal M$, we consider the fiberwise 
scalar product with a test section
$\zeta \in {\cal M}$.  As the supports are disjoint, it
is easy to see that

\[ \langle \xi_{total} | \zeta \rangle = \sum_{i = 1}^\infty 
\langle \xi_i|_{\neg \cup_{j < i} supp(\xi_j)} | 
\zeta \rangle \]

By Lemma \ref{res.lemma} each of the summands is
measurable, however,  it is quite easy to construct
examples in which the convergence of the sequence of
partial sums is not uniform a.e.

To see that the limit is, in fact, measurable, observe that the disjointness
of the supports of the summands implies that for any Borel set M, the sets
$N(f_i) \cap f_i^{-1}(M)$ are disjoint, where $f_i$ is the $i^{th}$
summand in the series of scalar products, and moreover, letting  
$f(x) = \langle \xi_{total}(x) | \zeta(x) 
\rangle_x$, that
$N(f) \cap f^{-1}(M)$
is their (disjoint)
union.  Note that the $N(f_i) \cap f_i^{-1}(M)$ are
measurable by Lemma \ref{res.lemma}, and thus, their
union is.

Having established that $\xi_{total}$ is a measurable
section of $\cal H$ with support equal to that of $\cal H$, observe that if
${\cal H}_x$ is non-zero, then the
singleton
$\{\xi_{total}(x)\}$ is a basis.  This suggests that we should
have an isomorphism from $\cal H$ to
${\Bbb C}|_A$ given by
the field of operators $\psi$ which maps $\xi_{total}(x)$ to
$1|_A$.  The difficulty is, that as it stands neither this nor
its inverse need be bounded fields of operators. 

To correct this, we replace
$\xi_{total}$ with $\xi_{normal} = \frac{1}{\|\xi_{total}(x)\|_x}\xi_{total}(x)$.  This is a measurable
section by Lemma \ref{measmult} below.  We can now give an
isomorphism from $\cal H$ to ${\Bbb C}|_A$ by the field of operators
$\phi$ which maps $\xi_{normal}(x)$ to $1|_A$.
 
The fields of operators $\phi$ and $\phi^{-1}$ in fact preserve the scalar
product in each fiber, and are thus plainly measurable and bounded.$\bullet$ \smallskip

Applied to any measurable field of Hilbert spaces $\cal H$ the construction 
just given shows that there is a 
measurable section
$\xi_{normal}$ with the same support as $\cal H$ and norm 1 in every
non-zero fiber..

This observation, together with Proposition \ref{split}
will show

\begin{thm} \label{partialin}
For every measurable field of Hilbert spaces $\cal H$,
there exists a partial measurable line bundle $\cal L$
which is a direct summand of $\cal H$ and has the
same support as $\cal H$.
\end{thm}

First observe

\begin{lemma} \label{measmult}
If $\xi$ is a measurable vector field for some measurable field $\cal H$ on 
$(X,S)$ and $\phi:X\rightarrow {\Bbb R}$ is a
measurable function, then $\phi \xi$, the section given
at $x$ by $\phi(x)\xi(x) \in {\cal H}_x$ is measurable.
\end{lemma}

\noindent{\bf proof:}  Let $\zeta \in {\cal M}_{\cal H}$,
then

\[ \langle \phi(x)\xi(x) | \zeta \rangle_x = 
\phi(x)\langle \xi(x) | \zeta \rangle_x \]

\noindent by linearity in each fiber.  But as products of 
measurable functions are measurable, the latter defines a measurable function 
on $X$, and by condition 2, $\phi \xi$ is a measurable vector field. $\bullet$ \smallskip

From this, together with the observation that sums, products and 
reciprocals of measurable real-valued functions are measurable, we see

\begin{propo} \label{gs} 
If the Gram-Schmidt process is applied (fiberwise) to a set of measurable 
vector fields, the result is a set of measurable vector fields.
\end{propo}

\begin{propo} \label{split}
If $\{\xi^{(i)}\}$ for $i = 1, . . ., n$ is a finite
set of measurable vector fields in $\cal H$ and
${\cal G}_x = span\{\xi^{(i)}(x)\}$ then 
${\cal G} = ({\cal G}_x, {\cal M}_{\cal H} \cap \prod_x
{\cal G})_x)$ is a measurable field of Hilbert spaces,
and moreover a direct summand of $\cal H$ with complementary summand given by
${\cal G}^\perp = ({\cal G}_x^\perp, {\cal M}_{\cal H} \cap \prod_x
{\cal G})_x^\perp)$.
\end{propo}

\noindent{\bf proof:}
By Proposition \ref{gs} it follows that the fiberwise orthogonal
projections onto $\cal G$ and ${\cal G}^\perp$ preserve
measurable sections.  Now, observe in both $\cal G$ and ${\cal G}^\perp$ 
the image of
the fundamental sequence of $\cal H$ under the orthogonal projections
may be taken as the fundamental sequence for the
subfield. Observe also that in either case, since the
projection is idempotent, it follows that the image of the
measurable sections is precisely ${\cal M}_{\cal H} \cap
\prod_x {\cal G}_x$ or ${\cal M}_{\cal H} \cap
\prod_x {\cal G}_x^\perp$ respectively.

We have thus established that conditions 1 and 3 hold in both
$\cal G$ and ${\cal G}^\perp$ , and that once condition 2 is shown, that the orthogonal
projections and inclusions are measurable fields of operators. 
The required equational condition for the direct sum decomposition follows
from fiberwise condition
in {\bf Hilb}.  The inclusions are plainly norm preserving, and it follows from
this, the equational condition, and orthogonality that the projections are norm
decreasing.

To establish condition 2, observe that any measurable 
vector field $\zeta$ decomposes as a sum of its 
projections onto measurable
vector fields $\zeta_\parallel \in {\cal G}$ and
$\zeta_\perp \in {\cal G}^\perp$.  By orthogonality, it
follows that for a section $\xi$ of $\cal G$ (resp. ${\cal G}^\perp$) the
scalar product $\langle \xi | \zeta \rangle$
is equal to $\langle \xi | \zeta_\parallel \rangle$
(resp. $\langle \xi | \zeta_\perp \rangle$).  From this and
condition 2 for $\cal H$, condition 2 for the subfields follows
immediately. $\bullet$ \newpage

\section{Bounded Invertible Additive Functors}

Since the primary motivation for this work is the construction of a suitable
setting for the representation theory of 2-groups, we begin considering
a family of functors sufficient for the construction of the 
representations
themselves.  In Section \ref{mf} we will construct the larger family of
functors which will be needed for the intertwiners in the theory
developed in \cite{CY.M2G}.  

\begin{defin}
An {\em invertible additive functor} is a functor $F$ between two
additive categories, which admits an inverse
functor (up to natural isomorphism) $G$ such that both $F$ and $G$ preserve
the addition of parallel maps, the zero object (up to canonical isomorphism)
and the biproducts (up to canonical isomorphism).
\end{defin}

\begin{defin}
A functor $F:{\cal C}\rightarrow {\cal D}$
between $C^*$-categories is {\em bounded} if there exists a 
constant $N > 0$ such that for all $f:X\rightarrow Y \in Arr({\cal C})$

\[ \|F(f)\| \leq N\|f\|. \]
\end{defin}

Observe that $*$-functors are bounded, being, in fact, norm decreasing
(cf. \cite{GLR}).

\begin{defin} \label{bnt}
A natural transformation $t$ between bounded functors is bounded
if $\|t\| = \sup_{X\in Ob({\cal C})} \|t_X\| < +\infty$.
\end{defin}

Our goal in this section is to characterize bounded invertible
additive functors between categories of measurable fields of Hilbert spaces
 in terms of the Borel space structures on the source and target.  
In particular, we will show that any bounded
invertible additive functor between
categories of the form ${\bf Meas}(X,S)$ is induced by an invertible
measurable transformation between the underlying Borel spaces.

The primary tool in showing this is partial measurable line bundles.
We begin by showing

\begin{thm} \label{pres.pmlb}
If $F:{\bf Meas}(X,S)\rightarrow {\bf Meas}(Y,T)$ is an
invertible additive functor with inverse $\Phi:{\bf Meas}(Y,T)
\rightarrow  {\bf Meas}(X,S)$,  and $\cal H$ is a partial
measurable line bundle on $X$, the $F({\cal H})$ is
a partial measurable line bundle on $Y$.
\end{thm}

\noindent{\bf proof:}  Suppose not.  Now by Theorem \ref{partialin}, 
$F({\cal H})$
admits a direct summand $\cal L$ which is a partial measurable line bundle 
with $supp({\cal L}) = supp(F({\cal H}))$.  Since $F(\cal H)$ is not itself a
partial measurable line bundle, ${\cal L}^\perp$ has non-empty (measurable)
support $A \subset supp({\cal H})$, and there exists a partial measurable
line bundle ${\cal K} \subset {\cal L}^\perp$ with support $A$.
Thus by Theorems \ref{split} and \ref{res.sum} we have a decomposition
of $F({\cal H})$ as

\[ F({\cal H}) \cong
{\cal L}|_A \oplus {\cal L}|_{\neg A}\oplus {\cal K}|_A
\oplus {\cal K}|_{\neg A} \oplus ({\cal L} + {\cal K})^\perp \]

\noindent
(Note the last summand is the orthogonal complement of $\cal K$ in
${\cal L}^\perp$.) But ${\cal L}|_A \cong {\cal K}|_A$, as they are both
partial measurable line bundles with support $A$. Now applying $\Phi$ to this
decomposition gives

\[ {\cal H} \cong \Phi(F({\cal H})) \cong \Phi({\cal L}|_A)
\oplus \Phi({\cal L}|_{\neg A})\oplus \Phi({\cal K}|_A)
\oplus \Phi({\cal K}|_{\neg A}) \oplus \Phi(({\cal L} + {\cal K})^\perp) ,\]

\noindent but since $\Phi({\cal L}|_A) \cong \Phi({\cal K}|_A)$
contains a partial measurable line bundle with the same support, the fibers of
$\cal H$ at points of $A$ have 
dimension greater than one, but this implies that
$supp(\Phi({\cal L}|_A) = \emptyset$,i.e. $\Phi({\cal L}|_A) = 0$.
But then, we have
${\cal L}|_A \cong F(\Phi({\cal L}|_A) \cong F(0) \cong 0$, the last by the
additivity of $F$.  But this is a contradiction, since 
$A \subset supp({\cal H})$ was
non-empty. $\bullet$ \smallskip
\smallskip

We also have

\begin{propo}If $F:{\bf Meas}(X,S)\rightarrow {\bf Meas}(Y,T)$ is an invertible
additive functor with inverse 
$\Phi:{\bf Meas}(Y,T)
\rightarrow  {\bf Meas}(X,S)$, and $\cal H$ and $\cal K$ are partial
measurable line bundles on $X$, then

\begin{enumerate}
\item $supp({\cal H}) \subset supp({\cal K})$ implies
$supp(F({\cal H})) \subset supp(F({\cal K}))$
\item $supp({\cal H}) \cap supp({\cal K}) = \emptyset$ implies
$supp(F({\cal H})) \cap supp(F({\cal K})) = \emptyset$
\end{enumerate}

\end{propo}

\noindent{\bf proof:}  For the first statement, use Theorem \ref{res.sum}
to decompose $\cal K$ into a direct summand isomorphic to $\cal H$ and a
direct summand supported on $\neg supp({\cal H})$.  Apply $F$, and observe that
the support of a direct summand is necessarily contained in the support of the
direct sum. For the second statement, suppose not, decompose
the direct sum of the images using Theorem \ref{res.sum} along the
non-empty measurable set $A = supp(F({\cal H})) \cap supp(F({\cal K}))$.
Taking images under $\Phi$ yields the same sort of contradiction as in
the proof of Theorem \ref{pres.pmlb}. $\bullet$ \smallskip

In the case of invertible additive functors, this result essentially
determines their structure:

\begin{cor} \label{induced.by.function}
If $F:{\bf Meas}(X,S)\rightarrow {\bf Meas}(Y,T)$ is an invertible additive
functor, then $F$ induces an isomorphism of $\sigma$-algebras
by $supp({\cal L}) \mapsto supp(F({\cal L}))$, which is induced
by pullback along an invertible measurable function 
$\widehat{F}:Y\rightarrow X$\end{cor}

Of course, it is easy to establish the converse:

\begin{propo} \label{equiv.from.function}
If $f:(Y,T)\rightarrow (X,S)$ is a measurable function with measurable
inverse, then it induces an additive equivalence of categories
between ${\bf Meas}(X,S)$ and ${\bf Meas}(Y,T)$ by pullback along
$f$ and its inverse.\end{propo}

The two previous results taken together suggest, and almost suffice to show

\begin{thm}
Any bounded additive isomorphism between categories of the form
${\bf Meas}(X,S)$ which has a bounded inverse
is naturally isomorphic to one induced by
an inverse pair of measurable functions.\end{thm}

\noindent{\bf proof:}  That the restriction of any such functor to the
full subcategory of partial measurable line bundles both induces such a
function on the underlying Borel spaces and is induced by one
follows easily from Corollary \ref{induced.by.function} and Proposition
\ref{equiv.from.function}.  

Now, as any equivalence of categories preserves all limits and
colimits which exist, it follows that the result holds if ${\bf Meas}(X,S)$
is replaced by the full subcategory of all measurable fields of Hilbert
spaces with a finite bound on the dimensions of their fibers, as these
are all finite limits (or colimits) of partial measurable line bundles.

Thus the only real work, and the only place where the boundedness 
hypotheses are needed, is in constructing the components of the
natural isomorphism at measurable fields of Hilbert spaces with
unbounded (possibly infinite) fiber dimensions.

Let $F:{\bf Meas}(X,S)\rightarrow {\bf Meas}(Y,T), F^{-1}:{\bf Meas}(Y,T)\rightarrow {\bf Meas}(X,S)$ be such an isomorphism.  The restriction of $F$ to partial measurable line bundles is induced by pullback
along an isomorphism of Borel spaces $\Phi$.  

We thus have another pair of inverse functors (easily seen to be $*$-functors
and thus {\em a fortiori} bounded) $\Phi^*$ and $\Phi^{-1*}$, and the restrictions
of $F$ and $\Phi^*$ to the full subcategories of partial measurable line bundles, or to
measurable fields of bounded dimension, coincide.

Let $\cal H$ be a measurable field of Hilbert spaces on $(X,S)$.  
As every
element of every fiber of $\cal H$ is contained in a partial measurable line bundle whose
inclusion may be taken to be norm preserving, and every
pair of elements is contained in a subfield of Hilbert spaces with fiber dimensions
less than or equal to two, it is easy to see that there is a linear isomorphism
from $F({\cal H})_y$ to $\Phi^*({\cal H})_y$ for all $y \in Y$, induced by the linear
isomorphisms between subfields with dimension bound two.  We need to see that this
family of linear isomorphisms, and the family of its inverses, in fact constitute
a bounded field of bounded operators on $(Y,T)$.  

Observe
that $\Phi^*$ is norm preserving--it is easy to see that pullbacks in general are
$*$-functors and thus norm decreasing, but being invertible, it must preserve norms.
So consider a vector $\zeta_y \in \Phi^*({\cal H})_y$.  It is a pullback of a vector
$\xi_x \in {\cal H}_x$ where $\Phi(y) = x$, which in turn lies in a section $\xi$
which generates a subfield $\iota_\xi: \langle \xi \rangle 
\rightarrow {\cal H}$, where the inclusion
$\iota_\xi$ may be taken to be norm preserving.  Applying $F$ and $\Phi^*$ to the inclusion
gives inclusions $F(\iota_\xi):F(\langle \xi \rangle)\rightarrow F({\cal H})$ and
$\Phi^*(\iota_\xi):\Phi^*(\langle \xi \rangle)\rightarrow \Phi^*({\cal H})$.

As $\Phi^*$ and $\iota_\xi$ are norm preserving in the relevant senses, $\Phi^*(\iota_\xi)$ 
is norm preserving.  But $F$ is bounded,
say with bound $N$, and $\iota_\xi$ is norm preserving and thus $F(\iota_\xi)$ is bounded
with bound $N$.  But the linear isomorphism constructed above carries $\zeta_y$ to the
image of $\zeta_y$ under $F(\iota_\xi)$ (recall that $F$ and $\Phi^*$ coincide on 
partial measurable line bundles).  Thus, the norm of $\zeta_y$ is dilated by at most a
factor of $N$.  As this applies to any vector in any fiber, we have shown that the family 
of linear isomorphisms is a bounded family of bounded operators.

To see that the family of inverses is a bounded family of bounded operators, observe
that an operator admits a bounded inverse exactly when its dilation of norms is bounded
away from zero, say by $1/N$,
and similarly that an invertible functor between $C^*$-categories has a bounded inverse if and only
if its dilations of norms are bounded away form $0$.  An identical argument using these
lower bounds then suffices to show that the inverse family is a bounded family of 
operators.
 
Naturality of the family of bounded fields of operators thus constructed follows
from the same considerations as linearity:  it can be checked elementwise, and all element
lie in partial measurable line bundles.
$\bullet$ \smallskip

This last result also allows us to characterize the natural transformations
between bounded invertible additive functors. 

First, to describe natural endomorphisms of such functors, observe
that by 1-composition with the inverse, it suffices to describe the
natural endomorphisms of the identity functor:

\begin{thm}
Any natural endomorphism of $Id:{\bf Meas}(X,S)\rightarrow {\bf Meas}(X,S)$
is given by multiplication by a bounded measurable function $\phi:X\rightarrow
{\Bbb C}$, and thus is bounded in the sense of Definition \ref{bnt}.
\end{thm}

\noindent{\bf proof:} By construction in the proof of the
previous theorem, any natural transformation
is determined by its components at partial measurable line bundles.  But by
Theorem \ref{res.sum}, these are determined by the component at the total
measurable line bundle ${\Bbb C}$.  Any map from ${\Bbb C}$ to itself is given
by multiplication by a bounded
measurable function $\phi$, and it is easy to see
that multiplication by any such function induces a natural endomorphism
of $Id$. $\bullet$ \smallskip

More generally the fact that the components on partial measurable line 
bundles determine natural transformations allows us to show

\begin{thm}
Natural transformations from $F:{\bf Meas}(X,S)\rightarrow {\bf Meas}(Y,T)$
to $G:{\bf Meas}(X,S)\rightarrow {\bf Meas}(Y,T)$, where both are invertible
additive functors are given by bounded measurable functions 
$\phi:E\rightarrow {\Bbb C}$, where $E$ is the equalizer of $\widehat{F}^{-1}$
and $\widehat{G}^{-1}$.
\end{thm}

\noindent{\bf proof:}Let ${\cal H} = {\Bbb C}|_A$ be a partial measurable 
line bundle
with support $A$.  Then $F({\cal H})$ (resp. $G({\cal H})$) is a
partial measurable line bundle on $Y$ with support $\widehat{F}^{-1}(A)$
(resp. $\widehat{G}^{-1}(A)$).  Thus maps between them are given by
multiplication by measurable functions on $\widehat{F}^{-1}(A) \cap
\widehat{G}^{-1}(A)$ (and zero on other fibers).

Thus, in particular a natural transformation must induce (and be induced
by) a family of measurable functions $\phi_A:\widehat{F}^{-1}(A)\cap
\widehat{G}^{- 1}(A)\rightarrow {\Bbb C}$.  However, not every such family
induces a natural transformation, as naturality squares impose consistency
conditions.  In particular, if we have a measurable set $B \subset A$,
the direct sum decomposition ${\Bbb C}|_A \cong {\Bbb C}|_B \oplus
{\Bbb C}|_{A\setminus B}$ provides naturality squares for the inclusions
and projections.  From these it follows that $\phi_A$ can be non-zero
only on $(\widehat{F}^{-1}(B) \cap \widehat{G}^{-1}(B)) \cup 
(\widehat{F}^{-1}(A\setminus B) \cap \widehat{G}^{-1}(A\setminus B))$.

Letting $B \subset A$ range over all containments of measurable subsets,
we find, that $\phi_A$ can only be non-zero on the 
points of $A$ which are actually in the equalizer of $\widehat{F}^{-1}$
and $\widehat{G}^{-1}$. 

It is easy to see that multiplication by any such function induces
a natural transformation.$\bullet$ \smallskip

%A particular class of additive reflections is of interest.  

%\begin{defin}
%A function between Borel spaces
%$f:(X,S)\rightarrow (Y,T)$ is {\em bimeasurable} if $A \in S$
%implies $f(A) \in T$, and $B \in T$ implies $f^{-1}(B)$ \in S$
%\end{defin} 

%Now, given a bimeasurable injection $i:(X,S)\rightarrow (Y,T)$,
%there exists a measurable retraction $i_*:{\bf Meas}(X,S)\rightarrow 
%{\bf Meas}(Y,T)$, $i^*:{\bf Meas}(Y,T)\rightarrow {\bf Meas}(X.S)$
%given by push 
%forward and pullback of measurable fields, respectively.

\section{Direct Integrals}

Classically the point of defining measurable fields of Hilbert 
spaces was to define
a measure theoretic analogue of direct sums for the purpose of decomposing
vonNeumann algebras and representations of non-compact groups.

In this section we investigate the functorial properties of this 
construction, and
extend it to families of objects in ${\bf Meas}(X,S)$.

Recall

\begin{defin}
Given a(n almost) measurable field $\cal H$ of Hilbert spaces on a Borel space
$(X,S)$ and a
measure $\mu$ on $(X,S)$, the {\em direct integral}

\[ \int^\oplus_X {\cal H}_x d\mu(x) \]

\noindent is the Hilbert space of all measurable sections 
$\xi \in {\cal M}_{\cal H}$ such that

\[ \|\xi \| = \{\int_x \|\xi(x) \|^2 d\mu(x)\}^{\frac{1}{2}} < +\infty \;\;\; (*)\]

\noindent modulo the identification of measurable sections which are equal 
$\mu$-a.e.
and equipped with the scalar product

\[ \langle \xi | \zeta \rangle = \int \langle \xi(x) | \zeta(x) \rangle d\mu(x)
 \]
 
We call a section $\xi$ satisfying $(*)$ an $L^2$ section of $\cal H$, and denote
its $\mu$-a.e. equivalence class by $\int^\oplus_X \xi(x) d\mu(x)$.

Similarly given a ($\mu$-essentially) bounded field of operators
$\alpha:{\cal H}\rightarrow {\cal K}$, the {\em direct integral}
$\int^\oplus_X \alpha(x) d\mu(x)$  is the
map which takes an element $\xi$ of $\int^\oplus_X {\cal H}_x d\mu(x)$  to the
element of $\int^\oplus_X {\cal K}_x d\mu(x)$ given at each point $x$ by
$\alpha(x)(\xi(x))$.

For a ($\mu$-essentially) bounded field of operators $\alpha(x)$, we denote the
map taking  $\int^\oplus_X \xi(x) d\mu(x)$ to the a.e.-equality class
of the section $x \mapsto \alpha(x)(\xi(x))$ by 
$\int^\oplus_X \alpha(x) d\mu(x)$.
\end{defin}

Now, it is easy to see that composition of bounded fields of 
operators is carried to
composition of operators, and that the identity field is carried to 
the identity
operator.  Thus, $\int^\oplus_X \, - \, d\mu(x)$ is a functor from 
${\bf Meas}(X,S)$ to ${\bf Hilb}$ for any measure $\mu$ on $(X,S)$.

It is also easy to see that two fields of operators which agree $\mu$-a.e. are
mapped to the same operator between the direct integrals.

Several categorical properties of this construction will be important

\begin{thm}
$\int^\oplus_X \;- \; d\mu(x)$ is a ${\Bbb C}$-linear additive functor from 
${\bf Meas}(X,S)$ to
${\bf Hilb} \cong {\bf Meas}(\{*\},\{ \{*\}, \emptyset \})$.  
Moreover, if $\mu$ is a probability
measure, the functor $\int^\oplus_X \;-\; d\mu(x)$ is bounded.
\end{thm}

\noindent {\bf proof:}

It is immediate by construction that this construction preserves identity
maps and carries the composition of bounded fields of operators to the
composition of operators.

For $\Bbb C$-linearity, first note that once it is observed that the sum
of two bounded fields of operators is a bounded field of operators, it is
immediate that it acts on $L^2$-sections to give the sum of the action
of its summands.  Likewise, multiplication of all fibers by a complex
scalar is a bounded field of operators, and induces the scalar multiplication
on the a.e.-equality classes of $L^2$ sections.

For additivity, note that the constant 0 field of Hilbert spaces is mapped
to the Hilbert space 0.  If we consider a biproduct in the category 
${\bf Meas}(X,S)$, the equational conditions required for the direct integral to
be a biproduct follow from functoriality and linearity.

Now let $\alpha$ be a bounded field of operators.  We then have

\[ \left\{ \int_X \|\alpha(x)(\xi(x))\|^2 d\mu(x) \right\}^\frac{1}{2} \leq
\sup_{x\in X} \|\alpha(x)\|_x 
\left\{ \int_X \|\xi(x)\|^2 d\mu(x) \right\}^\frac{1}{2} \]

But if $\mu$ is a probability measure we have

\[\left\{ \int_X \|\xi(x)\|^2 d\mu(x) \right\}^\frac{1}{2} \leq \sup_{x\in X} \|\xi(x)\|_x 
= \|\xi\| \]

\noindent and thus $\|\int_X^\oplus \alpha(x)(\xi(x)) d\mu(x) \| \leq \|\alpha\| \|\xi\|$. Therefore $\int_X^\oplus \; - \; d\mu(x)$ decreases the norm of arrows to which is applied,
and is thus a bounded functor.
$\bullet$ \smallskip
\smallskip

One result which will be quite useful, as it will allow us to reduce most 
proofs to the case of probability measures is the following:

\begin{thm}
If $\mu \ll \nu$ then there is a natural transformation from 
$\int^\oplus_X \;- \; d\nu(x)$ to $\int^\oplus_X \;- \; d\mu(x)$ induced
by multiplication by $\sqrt{\frac{d\mu}{d\nu}(x)}$, the square root of the
Radon-Nikodym derivative.

If $\mu \equiv \nu$, then this natural transformation is a natural isomorphism
with inverse induced by multiplication by  $\sqrt{\frac{d\nu}{d\mu}(x)}$.
\end{thm}

\noindent{\bf proof:} Recall that the Radon-Nikodym derivative
is always a measurable function, and that square-root of non-negative
functions preserves measurablity. Thus multiplication by 
$\sqrt{\frac{d\mu}{d\nu}(x)}$ takes measurable vector fields for $\cal H$ to 
measurable vector fields by Lemma \ref{measmult}. 

Thus the map on measurable fields of Hilbert spaces
takes a measurable vector field $\xi(x)$ to 
$\sqrt{\frac{d\mu}{d\nu}(x)}\xi(x)$.

But

\[ \int \|\sqrt{\frac{d\mu}{d\nu}(x)}\xi(x)\|_x^2 d\nu(x) = 
\int \|\xi(x)\|^2 \frac{d\mu}{d\nu}(x) d\nu(x) = \int \|\xi(x)\|_x^2 d\mu(x) \]

\noindent and thus the map induces a map on direct integrals as desired.

Naturality follows from the fact that the map is induced by a central
endomorphism of the measurable field of Hilbert spaces.

The second statement follows from the chain rule for Radon-Nikodym derivatives.
$\bullet$ \smallskip

Since any $\sigma$-finite measure is equivalent to a probability measure, we
have

\begin{cor} \label{equivtoprob}
Any direct integral functor $\int^\oplus_X (-)d\mu(x)$ is naturally 
equivalent to a direct integral functor $\int^\oplus_X (-) d\nu(x)$ for
a probability measure $\nu$ on $X$.
\end{cor}

We can generalize the direct integral construction to give rise to an important
family of functors between categories of the form ${\bf Meas}(X,S)$.
In fact, these functors will allow us to decompose many objects in categories
of measurable fields as direct integrals of simpler objects.

\begin{defin} \label{fiberedmeasure}
For a measurable function $\Phi:(X,S)\rightarrow (Y,T)$ between Borel spaces, 
a {\em $\Phi$-fibered measure on $X$} is a uniformly 
totally $\sigma$-finite conditional
measure distribution $\mu_y$, that is a $Y$-indexed family of measures
on $X$ such that

\begin{enumerate}
\item $\mu_y(X\setminus \Phi^{-1}(y)) = 0$
\item For all $A \in S$ the function $y \mapsto \mu_y(A)$ is measurable
\item There exist a sequence of measurable sets $A_n$ such that
$X = \cup_n A_a$ and
for all $y\in Y$ and all $n$ $\mu_y(A_n) < \infty$
\end{enumerate}
\end{defin}

We then have 

\begin{propo} If $f:X\rightarrow {\Bbb R}$ (or $\Bbb C$) is measurable and 
$\mu_y$ is a $\Phi$-fibered measure on $X$ for 
a measurable function $\Phi:X\rightarrow Y$, then
the function $y \mapsto \int f d\mu_y(x)$ is a measurable function on
$Y$
\end{propo}

\noindent{\bf proof:} It suffices to consider the case of $f$ real-valued
and non-negative.  Now, let $f_n$ be a sequence of simple functions
approximating $f$ from below.

The functions $y \mapsto \int f_n d\mu_y(x)$ are all measurable, being
real linear combinations of functions of the form $y \mapsto \mu_y(A)$
(for the $A$'s on which $f_n$ is constant).  But 
$[y \mapsto \int f d\mu_y(x)] = \limsup [y \mapsto \int f_n d\mu_y(x)]$, 
and thus is
measurable, as the $\limsup$ of a sequence of measurable functions is again
measurable.$\bullet$ \smallskip

This last result is useful to us because in generalizing direct
integrals to give functors between categories of the form ${\bf Meas}(X,S)$.
It will let us show that our constructions preserve measurable
sections.

We also have a fibered analog of our earlier reduction to probability
measures:

\begin{thm}If $\mu_y$ is a $\Phi$-fibered measure for 
a measurable $\Phi:(X,S)\rightarrow (Y,T)$, there exists
a $\Phi$-fibered measure $\nu_y$ such that for all $y$ $\mu_y \equiv
\nu_y$ and each $\nu_y$ is a probability measure on $X$.
\end{thm}

\noindent{\bf proof:}  The only catch in simply applying the standard
construction which shows that any $\sigma$-finite measure is equivalent
to a probability measure in each fiber separately is the need to
ensure that all the maps $y \mapsto \nu_y(A)$ for $A \in S$ will be
measurable.

Let $A_n$ be a sequence of disjoint measurable sets in $X$
such that $X = \cup_{n=1}^\infty A_n$ and all the $\mu_y(A_n)$'s are
finite.  Define $h:X\times Y\rightarrow {\Bbb R}$ by

\[ h(x,y) = \{2^n \mu_y(A_n)\}^{-1} \;\; x\in A_n. \]

We can then let $\nu_y(A) = \int_A h(x,y) d\mu_y(x) = \sum_{n=1}^\infty 
\{2^n \mu_y(A_n)\}^{-1} \mu_y(A_n \cap A)$.

But this is the $\limsup$ of its partial sums, which are measurable
since addition, multiplication,
multiplication by constants,
and reciprocation all preserve measurability, and the functions
$y\mapsto \mu_y(A_n)$ and $y\mapsto \mu_y(A_n \cap A)$ are both
measurable. Thus we are done. $\bullet$ \smallskip

We can then make

\begin{defin} \label{fibereddirectintegral}
Let $\Phi:(X,S)\rightarrow (Y,T)$ be a measurable function between 
Borel spaces.
Let $\mu_y(x)$ be a $\Phi$-fibered measure on $X$.

For a(n almost) measurable field of Hilbert spaces $\cal H$ on $(X,S)$,
let

\[ \int^\oplus_\Phi {\cal H}_x  d\mu_y(x) \]

\noindent 
denote the (almost) measurable field of Hilbert spaces with fiber at $y$
given by  $\int^\oplus_X {\cal H}_x d\mu_y(x)$, and
measurable sections given by the closure under condition 2 of the image
of the set of measurable sections of $\cal H$.
\end{defin}

Observe that this definition makes sense, since the fiberwise scalar
product on $\int^\oplus_\Phi {\cal H}_x  d\mu_y(x)$ is given by

\[ \left\langle \int^\oplus \xi(x) d\mu_y(x) |
\int^\oplus \zeta(x) d\mu_y(x) \right\rangle_y =
\int \langle \xi(x) | \zeta(x)
\rangle  d\mu_y(x) \]

\noindent which is the fiberwise integral of a measurable function
whenever $\xi$ and $\zeta$ are measurable sections of a(n almost) measurable
field of Hilbert spaces, and thus is measurable by the definition of
a $\Phi$-fibered measure.

In the measurable case the resulting field of Hilbert spaces admits a fundamental
sequence by virtue of the following two lemmas:

\begin{lemma}
If $\mu$ (resp. $\mu_y$) is a $\sigma$-finite measure (resp. $\Phi$-fibered measure) on
$X$,
then any measurable field of Hilbert spaces on $X$ admits a fundamental sequence of
$L^2$-sections with respect to $\mu$ (resp. of sections which are $L^2$-sections with respect to all of the $\mu_y$'s).  
\end{lemma}

\noindent {\bf proof:} Let $\{A_n\}_{n=1}^\infty$ be the sequence of measurable sets
exhausting $X$ of finite measure with respect to $\mu$ (resp. with respect to all of the
$\mu_y$).  Let $\{\xi_i\}_{i=1}^\infty$ be a fundamental sequence for $\cal H$.  Let 

\[ \zeta_{i,n,\alpha}(x) = \frac{\alpha}{\|\xi_i(x)\|} \xi_i|_{A_n}(x) \] 

\noindent for $i,n = 1, 2, . . .$ and $\alpha$ any algebraic number.  Plainly the 
$\zeta_{i,n,\alpha}(x)$ are dense in each ${\cal H}_x$.  But they are also $L^2$ sections
with respect to $\mu$ (resp. with respect to all $\mu_y$), since their norms are easily
seen to be given by

\[ \|\zeta_{i,n,\alpha}\| = \alpha \sqrt{\mu(A_n)} < +\infty \]

\noindent(resp. the same {\em mutatis mutandis} for each $\mu_y$).  Now, for technical
reasons we will replace this sequence with the sequence of all finite sums of its 
elements.
$\bullet$ \smallskip

\begin{lemma}
If $L^2(X,\mu)$ is separable (resp. $L^2(X,\mu_y)$ forms a measurable field of Hilbert
spaces on $Y$), then for any measurable field of Hilbert spaces ${\cal H}_x$ on $X$
which admits a fundamental sequence of $L^2$-sections with repect to $\mu$ (resp.
of sections which are $L^2$ with respect to all of the $\mu_y$),
the direct integral $\int^\oplus_X {\cal H}_x d\mu(x)$ is a separable Hilbert space
(resp. the almost measurable field of Hilbert spaces $\int^\oplus_X {\cal H}_x d\mu_y(x)$
on $Y$ is a measurable field of Hilbert spaces.)
\end{lemma}

\noindent{\bf proof:} First recall that the existence of a countable dense set in
a Hilbert space is equivalent to the existence of countable orthonormal basis.  A similar
argument establishes that the existence of a fundamental sequence in an almost 
measurable field of Hilbert spaces (i.e. the measurability of the field of Hilbert
spaces) is equivalent to the existence of a countable set of sections which are
orthogonal in each fiber, whose norms in all fibers is either 1 or 0. 
(Apply Gram-Schmidt fiberwise.)  Call such a countable set of sections a foundation.\footnote{
We want a word other than basis since it need not be a basis in each fiber, and may not
even be linearly independent globally.}

Consider $\{\phi_j\}$, an orthonormal basis for $L^2(X,\mu)$ (resp. a foundation for $L^2(X,\mu_y)$), and $\{\xi_i\}$, a foundation for $\cal H$.  We claim the set
$\{\int^\oplus_X \phi_j\xi_i d\mu \}$ (resp. $\{\int^\oplus_X \phi_j\xi_i d\mu_y \}$ )
is a total set in  $\int^\oplus_X {\cal H}_x d\mu(x)$ 
(resp. is total in  $\int^\oplus_X {\cal H}_x d\mu(x)$ for all $y$), and thus 
its algebraic span is dense (resp. dense in each
fiber).

Now, observe that 
$\| \phi_j(x)\xi_i(x) \|_x^2 = |\phi_j(x)|^2 \|xi_i(x)\|_x
\leq |\phi |^2 $.  Since $\phi_j$ is $L^2$ is it certainly measurable, and thus
$\phi_j \xi_i$ is a measurable sect
of its fiber-wise norm is majorized by the absolute square of the $L^2$ function
$\phi_j$.  It thus follows that $\phi_j\xi_i$ is an $L^2$-section.

It suffices to handle the case of a single measure space $(X,\mu)$.
Suppose $\int^\oplus_X \zeta d\mu$ is an element of $\{\int^\oplus_X \phi_j\xi_i d\mu \}$, and
that $\langle \int^\oplus_X \zeta d\mu | \int^\oplus_X \phi_j\xi_i d\mu \rangle = 0$ for
all $i$ and $j$.  But pulling out the $\phi_j$ by sequilinearity, we find that
 $\langle \zeta | \xi_i \rangle$ is orthogonal to all of the $\phi_j$'s, and thus is
zero in $L^2(X,\mu)$, since the $\phi_j$'s are total.
\footnote{As an aside, notice that this argument, together with 
the classical  proof of duality between $L^p$ spaces shows that 
when ever $\zeta$ and $\xi$ are $L^2$ sections of a measurable
field of Hilbert spaces, the function $x \mapsto \langle \zeta(x)
| \xi(x) \rangle_x$ is an $L^2$ function.}
Thus $\langle zeta | \xi_i \rangle$ is zero $\mu$-a.e., and thus 
$\zeta$ is zero $\mu$-a.e. since the $xi_i$'s are total in each fiber. Thus $\zeta$
is zero in the direct integral, and we are done. $\bullet$ \smallskip

At this point, we should point out why the preservation of measurability
under fiberwise integration is a non-vacuous condition, and provide some
useful examples of fibered measures.

We begin with an example of a family of measures on the inverse images
under a measurable map which is not a fibered measure:

\begin{exa}
Consider the second projection map $p_2:{\Bbb R}^2\rightarrow {\Bbb R}$.
Now, let $\beta:{\Bbb R}\rightarrow {\Bbb R}$ be a non-negative non-measurable
function.  Consider then the family of measures $\mu_y$ given by Borel
measure in the first co\"{o}rdinate $x$ on the subsets of
$[0,\beta(y)]\times \{y\} \subset p_2^{-1}(y)
\subset {\Bbb R}^2$, with sets lying in 
$p_2^{-1} \setminus [0,\beta(y)]\times \{y\}$
having measure zero.

The fiberwise integral of the (obviously measurable) constant function 1 on
${\Bbb R}^2$, which could arise, for instance, as the scalar product of
the constant section 1 in the constant field of Hilbert spaces ${\Bbb C}$
on ${\Bbb R}^2$ with itself,  is then the non-measurable function $\beta$.
\end{exa}

However, it is easy to give examples of fibered measures:

The following is essentially Fubini's Theorem for measurability:

\begin{exa} \label{measFub}
Let $(X,S)$ and $(Y,S)$ be any Borel spaces, and let $\mu$ be a measure 
on $X$, the family of measures given by 
$\mu_y(A) = \mu(p_1(A\cap (X\times{y})))$
is a $p_2$-fibered measure.  Condition 1 is true by construction, condition 2
is immediate since the functions $y \mapsto \mu_y(A)$ are all constant
functions, while for condition 3, if $A_n \subset X$ are a sequence of sets
which witnesses to the $\sigma$-finiteness of $\mu$, then $A_n\times X$ provide
the necessary sequence for the uniform $\sigma$-finiteness of $\mu_y$.
\end{exa}

Another example, which is the simplest example of a useful family we
will introduce later, is:

\begin{exa}
Let $(X,S)$ be any Borel space. Given $p_2:X\times X\rightarrow X$, 
the family of
measures $\mu_y$ on $X\times X$ given by

\[ \mu_y(A) = \left\{ \begin{array}{ll}
        1 & \mbox{if $(y,y) \in A$} \\
        0 & \mbox{otherwise}
        \end{array} \right. \]

\noindent
Again condition 1 is immediate, while any countable sequence of measurable sets
exhausting $X\times X$ will suffice for condition 3.  For condition 2, observe
that the projection functions are measurable, the diagonal $\Delta$ 
is measurable and
the function $y\mapsto \mu_y(A)$ in this case is simply the
characteristic function of $\Delta \cap A$.
\end{exa}

It is easy to see that this family is a fibered measure:  integration
with respect to it is simply restriction to the diagonal, which is a measurable
subset of $(X,S)$.

Similarly, given any section $s$ of the projection with a measurable image, the
family of
measures concentrated on $(s(y),y)$ giving each of these points 
measure 1 in its
fiber is a fibered measure.

$\Phi$-fibered measures are closed under multiplication by measurable functions
on the target, under addition, and suitable adaptations of limiting processes
which preserve measurable functions.

\section{Measurable Functors} \label{mf}

Having seen that ${\Bbb C}$-linear invertible additive functors
are induced by invertible
measurable functions between the underlying Borel space, we now turn to the
question of what class of functors including these are most appropriate
to consider when forming a 2-category of categories of (almost) measurable
fields of Hilbert spaces.

One obvious approach is to consider those ${\Bbb C}$-linear functors which
respect the norm structure in the sense that the map induced on
hom-sets is fiberwise continuous with respect to the operator norm.  

We will prefer a structural rather than an axiomatic approach--though we
conjecture that the class of functors just described is in fact the same
as that we are about to define.  

Given an
almost measurable field of Hilbert spaces ${\cal K}$ on $(X\times Y,S\star T)$,
and a $p_2$-fibered measure $\mu_y$ on $X \times Y$,  we may construct
a ${\Bbb C}$-linear additive functor $\Phi_{{\cal K},\mu_y}$:

Let 

\[ \Phi_{{\cal K},\mu_y }({\cal H})_y = 
	\int^\oplus_X {\cal H}_x \otimes {\cal K}_{<x,y>} d\mu_y(x) \]

\noindent with ${\cal M}_{\Phi({\cal H})}$ given as the closure
under condition 2 of the set 

\[ \left\{ \int^\oplus_X \eta(x) \otimes \kappa(x,y)
d\mu_y(x) | \eta \in {\cal M}_{\cal H}; \kappa \in {\cal M}_{\cal K} \right\} .
\]

Observe that there is a subtlety here:  not all pairs of sections 
$\eta, \kappa$ will define a section, only those which, for every $y\in Y$
give rise to $L^2$-sections of the tensor product.

\begin{defin}
A functor from ${\bf Meas}(X,S)$ to ${\bf Meas}(Y,T)$ is 
{\em measurable} if it is ${\Bbb C}$-linear 
equivalent to one of the form
$\Phi_{{\cal K},\mu_y}$. 
\end{defin}

We can then show the functors we considered earlier all belong to this
new class:

\begin{thm}
Any invertible additive functor between categories of the
form ${\bf Meas}(X,S)$ is a measurable functor.
\end{thm}

\noindent{\bf proof:} Let $F:{\bf Meas}(X,S)\rightarrow {\bf Meas}(Y,T)$
be an invertible additive functor. By Corollary \ref{induced.by.function}
it is induced by pullback along an invertible measurable function 
$\widehat{F}:(Y,T)\rightarrow (X,S)$.  If we then consider the constant
field of Hilbert spaces $\Bbb C$ on $X \times Y$ and the $p_2$-fibered measure

\[ \mu_y(A) = \left\{ \begin{array}{ll}
        1 & \mbox{if $(\widehat{F}(y),y) \in A$} \\
        0 & \mbox{otherwise}
        \end{array} \right. , \]

\noindent it is easy to see that the pullback functor is naturally 
isomorphic to the functor $\Phi_{{\Bbb C},\mu_y}$. $\bullet$ \smallskip

\begin{thm} If $\phi:(X,S)\rightarrow (Y,T)$ is a measurable function
between Borel spaces and $\mu_y$ is a $\phi$-fibered measure on $X$,
then the functor $\int^\oplus_\phi (-) d\mu_y(x)$ is a measurable functor.
\end{thm}

\noindent{\bf proof:}
Define a $p_2$-fibered measure on $X\times Y$ by 

\[ \tilde{\mu}_y(A) = \mu_y(\{ x | (x,\phi(x))\in A \}) \]

It is easy to see that $\int^\oplus_\phi (-) d\mu_y(x)$ is
naturally ismorphic to $\Phi_{{\Bbb C}, \tilde{\mu}_y}$, and thus
measurable. $\bullet$\smallskip

One vexing thing about measurable functors is the fact that it is
not immediate that the composition of measurable functors is measurable.

To show this, we will need to invoke Maharam's result \cite{Maharam}
on disintegration of measures, and several results relating tensor products
and direct integrals.

\begin{propo} \label{fieldspacetensor}
If $\cal H$ is a measurable field of Hilbert spaces on
$X$ and $\cal K$ is a separable Hilbert space, then there is
a measurable field of Hilbert spaces ${\cal H}\otimes {\cal K}$ 
with fiber at $x$ given by ${\cal H}_x \otimes {\cal K}$, with all 
algebraic-coefficient finite linear
combinations of tensor products of elements in the fundamental sequence
of $\cal H$ with elements of a countable dense set in $\cal K$ as 
fundamental sequence, and the closure under condition 2 of this 
sequence as ${\cal M}_{{\cal H}\otimes {\cal K}}$.
\end{propo}

\noindent{\bf proof:} By \cite{Takesaki} IV \S 8 Lemma 8.10, it suffices
to show that the proposed fundamental sequence is fiberwise dense and that
all of the pairwise scalar product functions are measurable.  Density
is clear from the construction of Hilbert space tensor products.

For the other condition observe that

\begin{eqnarray*}
\lefteqn{<\sum_{i,j} a_{i,j} \xi_j(x)\otimes \zeta_i | 
\sum_{k,l} b_{k,l} \phi_l(x)
\otimes \omega_k >_x =} \\
& &  \sum_{i,j,k,l} a_{i,j}\overline{b_{k,l}} 
<\zeta_i | \omega_k>_{\cal K} <\xi_j(x) | \phi_l(x)>_x , 
\end{eqnarray*}

\noindent is measurable as a function of $x$ since all of the $x
\mapsto  <\xi_j(x) | \phi_l(x)>_x$ are measurable, and linear combinations
of measurable functions are measurable (all other expressions occuring
in the last are constant in $x$). $\bullet$ \smallskip

Moreover we have

\begin{thm}
There is a canonical natural isomorphism

\[ \int^\oplus_X {\cal H}\otimes {\cal K}_x d\mu(x) \cong
\left( \int^\oplus_X {\cal H}_x d\mu(x) \right) \otimes {\cal K} \] 

\noindent where $K$ is any separable Hilbert space and ${\cal H}$ is 
any measurable field of Hilbert spaces on a Borel space $X$, and
${\cal H}\otimes {\cal K}$ is as in the previous proposition.
\end{thm}

\noindent{\bf proof:} 
By Corollary \ref{equivtoprob} we may assume without loss of 
generality that $\mu$ is a probability measure.

We proceed by first constructing a
canonical map from $\int^\oplus_X {\cal H}_x d\mu(x) \otimes {\cal K}$
to  $\int^\oplus_X {\cal H} \otimes {\cal K}_x d\mu(x)$, then showing
that it is an isomorphism and natural in both variables.

Now, since any element in  $\int^\oplus_x {\cal H}_x d\mu(x) \otimes {\cal K}$
is the limit of a Cauchy sequence of elements in the algebraic tensor
product, to specify a bounded linear operator from  
$\int^\oplus_x {\cal H}_x d\mu(x) \otimes {\cal K}$ to any other
Hilbert space, it suffices to specify its behavior on elements of the
form $\int^\oplus \xi(x) d\mu(x) \otimes \zeta$ for $\zeta \in {\cal K}$
and $\xi(x)$ a $L^2$-section of ${\cal H}$.

Rather obviously we wish to map  $\int^\oplus \xi(x) d\mu(x) \otimes \zeta$
to  $\int^\oplus \xi(x)\otimes \zeta d\mu(x)$.

To see that this in fact defines a bounded operator, we will need to
verify that

\begin{enumerate}
\item $\xi(x)\otimes \zeta$ is a measurable vector field in 
${\cal H}\otimes {\cal K}$,
\item $\{\int \|\xi(x)\otimes \zeta \|^2_x d\mu(x) \}^\frac{1}{2} < \infty$,
\item there exists and $M$ such that for all $\xi(x)$ and $\zeta$

 \[ \left\{\int \|\xi(x)\otimes \zeta \|^2_x d\mu(x) \right\}^\frac{1}{2} \leq
M \left\{\int \|\xi(x) \|^2_x d\mu(x) \right\}^\frac{1}{2}\|\zeta\|_{\cal K} \]

\noindent and
\item the image of each element is independent of the choice of
$\mu$-a.e. equality class representative $\xi(x)$.
\end{enumerate}

For the first consider an element of the fundamental sequence,
$\sum_{i,j} a_{i,j} \xi_i(x)\otimes \zeta_j$, and form the function of
$x$ given by scalar
product with $\xi(x)\otimes \zeta$. By sequilinearity this reduces
to a linear combination of the functions $x \mapsto <\xi_i(x) | \xi(x)>$,
and is thus measurable.

For the second and third, we compute

\begin{eqnarray*}
\left\{ \int \|\xi(x) \otimes \zeta \|^2 d\mu(x) \right\}^\frac{1}{2} & = &
        \left\{ \int \| \xi(x) \|^2_x \| \zeta \|^2_{\cal K} d\mu(x)
	\right\}^\frac{1}{2} \\
	& = & \|\zeta \| \left\{ \int \| \xi(x) \|^2 d\mu(x) 
\right\}^\frac{1}{2} 
\end{eqnarray*}

This is finite since the integral in the last right-hand side is the norm
of the $L^2$-section $\xi(x)$ of $\cal H$, that is the norm in the direct
integral $\int^\oplus {\cal H}_x d\mu(x)$.  Moreover, notice that the
right-hand side as a whole is the norm of the preimage $\int \xi(x) d\mu(x)
\otimes \zeta$, and thus $M = 1$ suffices.

For the independence of the choice of $\mu$-a.e. equality representative, 
observe that if $\xi = \xi^\prime$ $\mu$-a.e., then
$\xi \otimes \zeta = \xi^\prime \otimes \zeta$ $\mu$-a.e.

We have already shown more than that this map simply exists.  In showing
that it was bounded, we actually showed that it preserved the norm on a
dense set, and thus is an isometry onto its image.

It thus remains only to show that it is surjective.  But, by completeness,
it suffices to show that its image is dense in $\int^\oplus {\cal H}\otimes
{\cal K} d\mu(x)$.

Carrying out the proof of \cite{Takesaki} IV \S 8 Lemma 8.12 (simultaneous
fiberwise Gram-Schmidt orthonormalization) with the
fundamental sequence for ${\cal H}\otimes {\cal K}$ gives the same result
with one added feature:

\begin{lemma}
There is a sequence of measurable vectorfields $\psi_i$ in ${\cal H}
\otimes {\cal K}$ such that

\begin{enumerate}
\item $\{ \psi_i(x) | 1 \leq i \leq dim {\cal H}_x \}$ is an orthonormal
basis for ${\cal H}_x$,
\item for $i > dim {\cal H}_x$ $\psi_i(x) = 0$
\item $\psi_i = \sum_{k,l} a^i_{k,l} \xi_l \otimes \zeta_k$ for a finite
set of measurable functions $a^i_{k,l}(x)$, that is, $\psi_i$ lies in the
$Meas(X)$-linear span of the fundamental sequence of ${\cal H}\otimes 
{\cal K}$.
\end{enumerate}

\end{lemma}

Observe that the fundamental sequence of ${\cal H}\otimes {\cal K}$
lies in the image of the map.  It also follows from the bilinearity of
$\otimes$ and Lemma \ref{measmult} that all elements of the $Meas(X)$-linear
span are in the image of the map.

We can thus form a sequence of measurable subfields 

\[{\cal L}^n =
Span_{Meas(X)}\left\{ \psi_i | i = 1,. . .,n \right\}, \]

\noindent  and fields of bounded
operators $p_n$ and $p_n^\perp$ projecting onto these and their
orthogonal complements.  Moreover each of the ${\cal L}^n$ lies in the image of
the map.

Now, observe that $\psi = p_n(\psi) + p_n^\perp(\psi)$, and that
the summands are orthogonal in each fiber.

It thus follows that

\[ \|p_n(\psi)(x)\|^2 \leq \|\psi(x) \|^2 \]

\noindent and

\[ \|p_n^\perp (\psi)(x)\|^2 \leq \|\psi(x) \|^2 \]

By Proposition \ref{fieldspacetensor} $x \mapsto \| p_n(\psi)(x) \|^2$ 
is measurable, and since it is majorized by the integrable function
$x \mapsto \| \psi(x) \|^2$ it is integrable.  Thus 
$x \mapsto \| p_n^\perp(\psi)(x) \|^2$ is also integrable, being the 
difference of two integrable functions.

Each of $p_n(\psi)$ and $p_n^\perp(\psi)$
thus represents an element in $\int^\oplus {\cal H}_x \otimes {\cal K}
d\mu(x)$, and as noted above the $p_n(\psi)$ lie in the image of the
map.

It thus suffices to show that the sequence $\{ p_n(\psi) \}$ converge to
$\psi$, not merely pointwise in each fiber, but with respect to the norm on the
direct integral.

For any $\epsilon > 0$ let $E^\epsilon_n = \{ x \; | \;
\|\psi(x) - p_n(\psi)(x) \|_x^2 = \| p_n^\perp(\psi)(x) \|_x^2 < \epsilon \}$.

Observe that 

\begin{enumerate}
\item $m > n$ implies $E^\epsilon_m \supseteq E^\epsilon_n$;
\item all of the $E^\epsilon_n$ are measurable (being unions of sets
$N(f) \cap f^{-1}( (-\infty,\epsilon) )$ and $X\setminus (N(f) \cap f^{-1}
(\Bbb R))$, for $f(x) = \| p_n^\perp(\psi)(x) \|_x^2 $); and 
\item $\bigcup_{n=1}^\infty E^\epsilon_n = X$, since for all $x$ $\lim_{n\rightarrow
\infty} p_n(\psi)(x) = \psi(x)$.
\end{enumerate}

Thus we also have for any $\epsilon > 0$ that

$\int_{E^\epsilon_n} \| \psi(x) \|^2 d\mu(x)$ is an increasing sequence with
limit $M = \int_X \| \psi(x) \|^2 d\mu(x)$.

Now, fix $\epsilon > 0$. By the foregoing discussion, we can choose an 
$N$ such that

\[ M - \int_{E^\frac{\epsilon}{2}_N} \| \psi(x) \|^2 d\mu(x) 
< \frac{\epsilon}{2} . \]

\noindent
For any $n \geq N$, we then have

\begin{eqnarray*}
\| \psi - p_n(\psi) \|^2_{{\cal H}\otimes {\cal K}} & = & \int_X \| \psi(x) - p_n(\psi)(x) \|_x^2 d\mu(x) \\
       & = &  \int_{E^\frac{\epsilon}{2}_N} \| \psi(x) - p_n(\psi)(x) \|_x^2 d\mu(x) + \\
       &  & \hspace{1cm} \int_{X\setminus E^\frac{\epsilon}{2}_N} \| \psi(x) - p_n(\psi)(x) \|_x^2 d\mu(x) .
\end{eqnarray*}

\noindent But

\begin{eqnarray*}
\int_{E^\frac{\epsilon}{2}_N} 
       \| \psi(x) - p_n(\psi)(x) \|_x^2 d\mu(x) & < & 
       \int_{E^\frac{\epsilon}{2}_N} \frac{\epsilon}{2} d\mu(x) \\
       & \leq &
       \frac{\epsilon}{2} .
\end{eqnarray*} 

\noindent the first inequality by the construction of $E^\frac{\epsilon}{2}_N$,
and the fact that $E^\frac{\epsilon}{2}_N \subset E^\frac{\epsilon}{2}_n$,
the second because we are integrating over a probability space, while

\begin{eqnarray*}
 \int_{X\setminus E^\frac{\epsilon}{2}_N} 
	\| \psi(x) - p_n(\psi)(x) \|_x^2 d\mu(x) & = &
	 \int_{X\setminus E^\frac{\epsilon}{2}_N} 
	\| p_n^\perp(\psi)(x) \|_x^2 d\mu(x) \\  
	& \leq &
	 \int_{X\setminus E^\frac{\epsilon}{2}_N} 
	\| \psi(x) \|_x^2 d\mu(x) \\ 
	& = & M - \int_{E^\frac{\epsilon}{2}_N} \| \psi(x) \|^2 d\mu(x) \\ 
	& < & \frac{\epsilon}{2} .
\end{eqnarray*}

Thus the image is dense, and the map is an isomorphism.

Naturality in both variables follows easily by chasing the images of 
elements of the form $\int^\oplus \xi(x) d\mu(x) \otimes \zeta$. $\bullet$ \smallskip

\begin{thm} The composition of two measurable functors is a measurable
functor.
\end{thm}

\noindent{\bf proof:} Consider the functors $F = \Phi_{{\cal F}, \mu_y}$
and $G = \Phi_{{\cal G}, \nu_z}$, for $\cal F$ (resp. $\cal G$) 
a measurable field of Hilbert spaces on $X\times Y$ (resp. $Y\times Z$)
and $\mu_y$ (resp. $\nu_z$) a $p_2$-fibered measure on $X\times Y$
(resp. $Y\times Z$).  Without loss of generality we may assume that all
of the $\mu_y$ and $\nu_z$ are probability measures.

We then have

\[ G(F({\cal H}))_z = \int^\oplus_Y \left\{ \int^\oplus_X {\cal H}_x
\otimes {\cal F}_{(x,y)} d\mu_y(x) \right\} \otimes {\cal G}_{(y,z)} d\nu_z(y)
 \]

By the functoriality of the outer direct integral and the previous
theorem, this is then naturally isomorphic to

\[  \int^\oplus_Y \int^\oplus_X {\cal H}_x
\otimes {\cal F}_{(x,y)} \otimes {\cal G}_{(y,z)} d\mu_y(x) d\nu_z(y)
 \]

Now, for each $z \in Z$ $\int \mu_y d\nu_z(y)$ is a probability measure
on $X \times Y$.  Define a $z$-indexed
family of measures on $X$ by $\lambda_z(A) = \int \mu_y d\nu_z(y) (A\times Y)$.
By abuse of notation, we also denote by $\lambda_z$ the family of
measures on $X\times Z$ given by $\lambda_z(B) = \lambda_z(p_1(B\cap (X\times
\{ z \})))$, 
where $\lambda_z$ on the right-hand side is the measure just defined.

\begin{lemma}
The family of measures $\lambda_z$ on $X\times Z$ is a $p_2$-fibered measure.
\end{lemma}

\noindent{\bf proof:} By construction the first condition of Definition
\ref{fiberedmeasure} is satisfied.  The second follows from the corresponding 
condition for $\mu_y$ and $\nu_z$, while the third is immediate once it
is observed that the $\lambda_z$ are all probability measures.
$\bullet$ \smallskip

Observe also that by construction $\int \mu_y d\nu_z(y)$ satisfies the 
hypotheses of the following theorem of Maharam \cite{Maharam} (cf. also
\cite{GrMa}) with respect to the projection $p_1$ onto $X$ and
each of the $\lambda_z$:

\begin{thm} {\rm (Maharam)} 
Let $L$ be a Lusin space and $S$ a non-empty Suslin space, $p$ a measurable
map from $L$ to $S$ each equipped with the usual Borel structure.
A $\sigma$-finite measure $\mu$ on the Borel sets
of $L$ has a uniformly $\sigma$-finite disintegration
with respect to $p$ and $\nu$ a measure on $S$ if and only
for all measurable $B$ in $S$ $\nu(B) = 0$ implies $\mu(p^{-1}(B)) = 0$
\end{thm}

We thus have a $Z$-indexed family of $p_1$-fibered measures $\kappa_{x,z}$
such that

\[ \int \kappa_{x,z} d\lambda_z(x) = \int \mu_y d\nu_z(y) \]

Thus the object

\[  \int^\oplus_Y \int^\oplus_X {\cal H}_x
\otimes {\cal F}_{(x,y)} \otimes {\cal G}_{(y,z)} d\mu_y(x) d\nu_z(y)
 \]

\noindent can also be described as

\[  \int^\oplus_Y \int^\oplus_X {\cal H}_x
\otimes {\cal F}_{(x,y)} \otimes {\cal G}_{(y,z)} d\kappa_{x,z}(y) 
d\lambda_z(x)
 \]

\noindent Thus, applying the previous theorem and the functoriality
of the outer direct integral, we see that this is naturally
isomorphic to

\[  \int^\oplus_Y {\cal H}_x
\otimes \left[ \int^\oplus_X  {\cal F}_{(x,y)} \otimes {\cal G}_{(y,z)} 
d\kappa_{x,z} \right] d\lambda_z(x) .
 \]

Thus the composition is measurable, being naturally isomorphic to
the functor $\Phi_{\cal K},\lambda_z$ where

\[ {\cal K}_{(x,z)} = \int^\oplus_Y {\cal F}_{(x,y)}\otimes {\cal G}_{(y,z)}
d\kappa_{x,z} = \int^\oplus_{p_{1,3}} p^*_{1,2}{\cal F}_{(x,y)} \otimes
p^*_{2,3}{\cal G}_{(y,z)} d\kappa_{x,z} \]

is a measurable field by the construction of Definition 
\ref{fibereddirectintegral}.
$\bullet$\smallskip

In the same way, we can isolate an interesting class of natural 
transformations between measurable functors:

Consider two parallel measurable functors 

\[ \Phi_{{\cal F},\mu_y},
\Phi_{{\cal G},\nu_y}:{\bf Meas}(X,S)\rightarrow {\bf Meas}(Y,T). \]

Now if the fibered measures $\mu_y$ and $\nu_y$ are equal, it is
clear that any measurable field of operators on $X\times Y$ from
$\phi:{\cal F}\rightarrow {\cal G}$, 
which is essentially bounded with respect to $\nu_y$ on each
$X\times \{y\}$ will induce a natural tranformation with 
components given by

\[ \Phi_{\phi,\nu_y,{\cal H}} = \int^\oplus Id_{{\cal H}_x}\otimes \phi_{(x,y)}
d\nu_y(x). \]

If we restrict our attention to 
measures with $\mu_y$ is absolutely continuous with respect to $\nu_y$ for
every $y\in Y$, we can also
use any essentially bounded field of operators to induce a natural 
transformation, but only after we normalize by multiplying by
$\sqrt{\frac{d\mu_y}{d\nu_y}}$, the square root of the fiberwise
Radon-Nikodym derivative of $\mu_y$ with respect to $\nu_y$.

Even for pairs of totally $\sigma$-finite measures, for which
the measure defining the source is not absolutely continuous
with respect to the measure defining the target, we
can apply the Lebesgue decomposition theorem to construct a natural
transformation induced by any $\nu_y$ essentially bounded field
of operators: Decompose $\mu_y$ as $\tilde{\mu}_y + \hat{\mu}_y$,
where $\tilde{\mu}_y$ is dominated by $\mu_y$ and is 
absolutely continuous with respect to
$\nu_y$ and $\hat{\mu}_y$ and $\nu_y$ are singular.  Now observe
that for any section $\zeta_{(x,y)}$ of ${\cal H}\otimes {\cal K}$,
the condition that $\{\int^\oplus \|\zeta_{x,y}\|^2 d\mu_y(x)\}^\frac{1}{2}$
be finite implies the same condition with $\mu_y$ replaced with
$\tilde{\mu}_y$.  We can thus apply the construction of the previous
paragraph to use any $\nu_y$-essentially bounded field of operators
to induce a natural transformation.  (Note:  there is in general an
non-trivial kernel (in the algebraic sense) 
in passing from the direct integral with respect to
$\mu_y$ to that with respect to $\tilde{\mu}_y$.)

\begin{defin}
A {\em measurable natural transformation} between two measurable functors

\[
\Phi_{{\cal F},\mu_y},
\Phi_{{\cal G},\nu_y}:{\bf Meas}(X,S)\rightarrow {\bf Meas}(Y,T) \]

\noindent is a natural transformation with component at $\cal H$ given by

\[ \zeta_y \mapsto \int^\oplus 
\sqrt{\frac{d\tilde{\mu}_y}{d\nu_y}}
Id_{{\cal H}_x}\otimes B_{(x,y)}(\zeta_{x,y})
d\nu_y(x) \]

\noindent where $\zeta_y = \int^\oplus \zeta_{x,y} d\mu_y(x)$, for 
some field of operators $B:{\cal F}\rightarrow {\cal G}$ which is
$\nu_y$-essentially bounded for all $y$.
\end{defin}

It is easy to see that identity natural transformations are measurable,
and that both compositions of measurable natural transformations are
again measurable (this latter needing the chain rule for Radon-Nikodym
derivatives). 

We denote the bicategory whose objects are the categories 
${\bf Meas}(X,S)$, 1-arrows are all measurable functors, and 2-arrows
are all measurable natural transformations by {\bf Meas}.

\section{Tensor Products}

We have already considered a number of constructions involving tensor 
products, all of which have been quite well-behaved.  The following
result is thus not surprising:

\begin{thm}
For any Borel space $(X,S)$ the category ${\bf Meas}(X,S)$ is a
monoidal category when equipped with

\[ [{\cal H}\otimes {\cal K}]_x = {\cal H}_x \otimes {\cal K}_x \]

\noindent with fundamental sequence given by all algebraic linear combinations
of elements of the form $\eta_i \otimes \kappa_j$, where $\{\eta_i \}$
and $\{\kappa_j \}$ are fundamental sequences for $\cal H$ and
$\cal K$ respectively, the measurable sections are the closure
of this fundamental sequence under condition 2, and the total measurable
line bundle as $I$.  Structure maps are given fiberwise by the corresponding
structure maps for Hilbert-space tensor product.
\end{thm}

\noindent{\bf proof:}  The condition from \cite{Takesaki} IV \S 8 Lemma 8.10
follows from the same argument as in Proposition \ref{fieldspacetensor}.
Coherence follows from the coherence for the structure maps in each fiber.
$\bullet$\smallskip

What is perhaps a little more surprising is that these tensor products
and the cartesian product of Borel spaces induce a monoidal bicategory
structure on the 2-category {\bf Meas}:

\begin{thm}
{\bf Meas} is a mononidal bicategory when equipped with the
monoidal bifunctor $\odot$ given on objects by
${\bf Meas}(X) \odot {\bf Meas}(Y) = {\bf Meas}(X\times Y)$,
and induced on 1- and 2-arrows by the functors
$\odot = p_X^*(p_1)\otimes p_Y^*(p_2) :{\bf Meas}(X)\times 
{\bf Meas}\rightarrow {\bf Meas}(X\times Y)$, 
where $p_X$ and $p_Y$ are the projection maps from
$X\times Y$ onto its factors and $p_1$ and $p_2$ are the projection functors 
from ${\bf Meas}(X)\times {\bf Meas}(Y)$ onto its factors, and the object
${\bf 1}$ as identity.
\end{thm}

\noindent{\bf sketch of proof:} The structural 1-arrows are induced by the
corresponding structural arrows for cartesian product of Borel spaces and
tensor product of fields of Hilbert spaces.  The structural 2-arrows in
turn are all identities either by the coherence for the two products inducing
the structural 1-arrows, or by virtue of the functoriality in each
variable of the operation inducing $\odot$. $\bullet$

\section{Direct Integrals in ${\bf Meas}(X)$ }

The construction of measurable functors given above may be used to
construct direct integrals of $Y$-indexed families of objects (for
a Borel space $Y$) in any
${\bf Meas}(X)$.

\begin{defin}
A {\em measurable field of ${\bf Meas}(X)$-objects} on a Borel space
$Y$ is a measurable field of Hilbert spaces on $X\times Y$.  Similarly
a (bounded) field of ${\bf Meas}(X)$-arrows on $Y$ is a (bounded) 
field of operators on $X\times Y$.
\end{defin}

We can then define a direct integral of such a measurable field by

\begin{defin}
Given a measurable field $\cal K$ of ${\bf Meas}(X)$-objects on $Y$, 
and a measure
$\nu$ on $Y$, the direct integral 

\[ \int^\oplus {\cal K}_{<x,y>} d\nu(y) \]

\noindent is the image of the total measurable line bundle $\Bbb C$ on $Y$ under
the measurable functor $\Phi_{{\cal K},\nu}$, where $\nu$ is interpreted
as the $p_1$ fibered measure for which

\[ \nu_x(A) = \nu(p_2(A \cap p_1^{-1}(x)) . \]
\end{defin}

Observe that here we have switched the role of first and second projection, 
however the fact that this defines a measurable field follows from the
work done in the section on measurable functors.

This definition, together with the fibered measure of Example \ref{measFub}
allows us to give a version of Fubini's Theorem for direct integrals.

\begin{thm}
If $\mu$ (resp. $\nu$) is a measure on the Borel space $X$ (resp. $Y$),
then for any measurable field of Hilbert spaces $\cal H$ on $X \times Y$

\[ \int^\oplus_{X\times Y} {\cal H}_{(x,y)} d\mu(x)\times d\nu(y) \cong
\int^\oplus_X \left\{ \int^\oplus_Y {\cal H}_{(x,y)} d\nu(y) \right\} d\mu(x) .
\]
\end{thm}

\noindent{\bf proof:} The result follows from the classical Fubini's Theorem
applied to the integrals in the definining conditions for measurability
and $L^2$-ness of sections. $\bullet$

\section{Measurable Categories and Measurable Bicategories}

As observed in the introduction, our purpose in examining in detail the
structure of categories of measurable fields of Hilbert spaces and of
organizing them into a (monoidal) bicategory was to provide a setting
for a representation theory for categorical groups.

That theory is developed in \cite{CY.M2G} and \cite{CS.2G}. It will
turn out to be the first example of what should be a general theory
of measurable bicategories.  Our purpose in this final section is
to suggest the outline of general theories of measurable categories
over a measure space $X$ and of
measurable bicategories which are analogues of Tannakian categories,
but with VECT replaced by ${\bf Meas}(X)$ and ${\bf Meas}$ respectively.

\begin{defin}
{\em A measurable category over $X$} is a monoidal category $\cal C$ equipped
with monoidal functors $U:{\cal C}\rightarrow {\bf Meas}(X)$ and 
$T:{\bf Meas}(X)\rightarrow {\cal C}$ such that there is
a natural isomorphism $U(T) \cong Id_{{\bf 
Meas}(X)}$.  Objects of $\cal C$ isomorphic to object of the form $T({\cal H})$
are called {\em trivial objects}.
\end{defin}

\begin{exa}
${\bf Meas}(X)$ with both structural functors being the identity functor.
\end{exa}

\begin{exa}
Fix a Lie group $G$, let ${\bf Rep}(G)/X$ be the category whose objects
are measurable fields of Hilbert spaces on $X$ such that each fiber is
equipped with a unitary representation of $G$, and moreover for each
$g \in G$ $g$ acts by a bounded field of bounded operators.  Then
the forgetful functor to underlying measurable fields and the inclusion
on measurable fields with a trivial $G$-action on each fiber make
${\bf Rep}(G)/X$ into a measurable category over $X$.
\end{exa}

Of more interest is the general setting in which the motivating construction
fits:

\begin{defin}
{\em A measurable bicategory $\cal C$} is a monoidal bicategory $\cal C$
equipped with monoidal bifunctors $U:{\cal C}\rightarrow {\bf Meas}$
and $T:{\bf Meas}\rightarrow {\cal C}$ such that $U(T)$ is naturally
isomorphic to $Id_{\bf Meas}$
\end{defin}

Of course {\bf Meas} itself gives a tautological example.

The subject of \cite{CY.M2G} gives others:

\begin{exa} Let $\cal G$ be a categorical group.  Regarding $\cal G$
as a bicategory with one object, the functor bicategory ${\bf Meas}^{\cal G}$
(with bifunctors as objects, pseudonatural transformations as 1-arrows,
and modifications as 2-arrows) is a measurable bicategory (with monoidal
structure induced by $\odot$ in an obvious way).
\end{exa}

The subcategories considered in \cite{CY.M2G} provide more examples.

Somewhat curiously, letting the underlying Borel space vary, measurable
categories themselves can be organized into a measurable bicategory:

\begin{exa}{\bf MeasCat} is a measurable bicategory whose objects
are all measurable categories; whose 1-arrows from $\cal C$ to $\cal D$
are pairs of monoidal functors $F:{\cal C}\rightarrow {\cal D}$ and
$\Phi:{\bf Meas}(X)\rightarrow {\bf Meas}(Y)$ with $\Phi$ measurable
and such that both the
obvious square formed by $F$, $\Phi$ and the underlying functors
and the obvious square formed by $F$, $\Phi$ and the inclusion of
trivial object functors commute; and whose 2-arrows from $(F,\Phi)$
to $(G,\Gamma)$ are pairs of 2-arrows $(h:F\Rightarrow G, \eta:\Phi\Rightarrow \Gamma)$
such that the two obvious pillows commute.

The underlying functor simply assigns ${\bf Meas}(X)$ to 
the measurable category $({\cal C},{\bf Meas}(X),U,T)$, while the
inclusion of trivial objects assigns to ${\bf Meas}(X)$, the 
tautological $({\bf Meas}(X),{\bf Meas}(X),Id,Id)$.

The monoidal structure on {\bf MeasCat} is given by letting
${\cal C}\odot {\cal D}$, for ${\cal C}$ (resp. $\cal D$) a measurable
category over $X$ (resp. $Y$), be the measured category over
$X\times Y$ with objects given by a pair of a $Y$-indexed family of $\cal C$
objects whose underlying $Y$-indexed family of fields of Hilbert
spaces on $X$ forms a measurable field of Hilbert spaces on $X\times Y$
and an $X$-indexed family of $\cal D$ objects whose underlying $X$-indexed
family of fields of Hilbert spaces on $Y$ forms a measurable field of Hilbert
spaces on $X\times Y$.  The underlying field of Hilbert spaces for
an object $(C_y, D_x)$ is the tensor product of the two underlying fields
of the entries.  The arrows are generated by the family similarly 
defined, and by formally adjoined isomorphisms between 
$(C_y\otimes T_{\cal C}^Y({\cal H}_<x,y>), D_x)$ 
and $(C_y, T_{\cal D}^X({\cal H}_<x,y> \otimes D_x)$, 
where $T_{\cal C}^Y$ (resp. $T_{\cal D}^X$) is the trivial functor
for $\cal C$ (resp. $\cal D$) applied in each of fiber over $Y$ (resp. $X$).
\end{exa}

\appendix

\section{Elementary results on ${\bf Hilb}$}

One problem in approaching this work is a paucity of references on the 
categorical structure of the category of Hilbert spaces and bounded 
operators.  The results in this appendix are elementary and, in
non-categorical guise, classical.  We include them for completeness
of exposition, but in an appendix so as not to interrupt the flow of 
the new results.

\begin{defin}
{\em The category of separable Hilbert spaces} ${\bf Hilb}$ has as objects
all separable Hilbert spaces and as arrows all bounded operators.  Source,
target, identities and composition are obvious.
\end{defin}

Many facts about ${\bf Hilb}$ are set forth in \cite{GLR}.

Of these, the most important for us is the fact that ${\bf Hilb}$ is a
$C^*$-category.  It is thus {\em a fortiori} equivalent to its 
opposite category by an equivlence which is the identity on objects,
adjoint-operator being the functor in the equivalence in either 
direction.

Ghez, Lima and Roberts \cite{GLR}, do not, however address some of
the elementary category-theoretic properties we will use.  We
summarize these in 

\begin{propo}
${\bf Hilb}$ is a ${\Bbb C}$-linear additive category with all finite limits
and colimits.
\end{propo}

\noindent{\bf proof:} ${\Bbb C}$-linearity is obvious.  To see
that ${\bf Hilb}$ is additive, observe that the $0$-dimensional vectorspace
is a Hilbert space and is a zero (initial and terminal) object for 
${\bf Hilb}$, likewise the vector space direct sum of two Hilbert spaces,
equipped with the scalar product $\langle (\xi, \phi) | (\zeta, \psi) \rangle
= \langle \xi, \zeta \rangle + \langle \phi, \psi \rangle$ is a Hilbert 
space, and all of the
structual maps for it as a vector-space biproduct are bounded operators.
Thus it is a biproduct in ${\bf Hilb}$. 

 By self-duality, it suffices to show that ${\bf Hilb}$
has all finite limits, and by standard results, it suffices to
show that it has kernels and binary products. Kernels are easy:  the
vector-space kernel will again be a Hilbert space
with the scalar product inherited from the source of the operator.  (Note:
since bounded operators are continuous, the kernel will be closed, and thus
complete.)  Finite products follow from biproducts.  $\bullet$ \smallskip

\begin{flushleft}
\bibliography{MeasCatand2Groups}
\end{flushleft}

\end{document}